\documentclass[usebiblatex=false]{kulifac}
\usepackage[mode=pub]{arxivix}

\DeclareUnicodeCharacter{2212}{\textendash}

\newcommand{\policy}{\pi}

\graphicspath{{assets/}}

\usetikzlibrary{patterns,hobby}

\begin{document}
\begin{frontmatter}

\title{Policy iteration using Q-functions: \\Linear dynamics with multiplicative noise.\thanksref{footnoteinfo}} 

\thanks[footnoteinfo]{This work is supported by: the Research Foundation Flanders (FWO) PhD grant 11E5520N and research projects G081222N, 
G033822N, G0A0920N; European Union's Horizon 2020 research and innovation programme under the Marie Skłodowska-Curie 
grant agreement No. 953348; Research Council KU Leuven C1 project No. C14/18/068; Fonds de la Recherche Scientifique
- FNRS and the Fonds Wetenschappelijk Onderzoek - Vlaanderen under EOS project no 30468160 (SeLMA).}

\author[First]{Peter Coppens} 
\author[First]{Panagiotis Patrinos} 

\address[First]{Department of Electrical
Engineering (ESAT-STADIUS), KU Leuven, Kasteelpark Arenberg
10, 3001 Leuven, Belgium.
{Email: \tt\normalsize peter.coppens@kuleuven.be, panos.patrinos@kuleuven.be}}

\begin{abstract}                % Abstract of not more than 250 words.
This paper presents a novel model-free and fully data-driven policy iteration scheme for quadratic regulation of linear dynamics 
with state- and input-multiplicative noise. The implementation is similar to the \emph{least-squares temporal difference} 
scheme for Markov decision processes, estimating Q-functions by solving a least-squares problem with instrumental variables. 
The scheme is compared with a model-based system identification scheme and natural policy gradient through numerical experiments. 
\end{abstract}


\end{frontmatter}
%===============================================================================

\section{Introduction}
The success story of AlphaZero in \cite{Silver2017a} and similar reinforcement learning strategies in a variety
of applications has reinvigorated the interest of the control community in learning schemes. Many of these 
fall under the general framework of dynamic programming, policy iteration, etc. (see \cite{Bertsekas2022, Bertsekas2022a} for a good overview). 
Hence understanding the performance of such schemes in simple settings can generalize to more complex problems. 
One setting under investigation has been \emph{Linear Quadratic Regulation (LQR)}. For example \cite{Recht2018,Mania2019,Dean2019} consider \emph{model-based} 
approaches where system identification is performed separately from control synthesis. 
In \cite{Bradtke1994,Lewis2009,Bu2019b} meanwhile \emph{model-free} approaches -- more akin to usual reinforcement learning -- like policy iteration and policy gradient --
were investigated. 

A setting with similar potential is quadratic regulation of linear dynamics with multiplicative noise \citep{Wonham1967}. Here analytical 
solutions are still available, yet several challenges absent from classical LQR present themselves. For example, the optimal 
controller depends on the first and second moment of the noise -- unlike in LQR with additive noise -- and the separation principle does not hold.
Moreover, as we will see, data-driven policy evaluation is not possible through normal least-squares as in \citep{Bradtke1994} 
and instead requires instrumental variables as is the case with Markov decision processes \citep{Bradtke1996}. Multiplicative
noise also arises naturally in aerospace and vehicle control applications \citep[\S1.9]{Damm2004}, biological applications \citep{Todorov2005,Mohler1980b} and
communication channels \citep{Wang2002} among others. Finally, the inclusion of multiplicative noise also induces robustness 
against parametric uncertainty \citep{Coppens2022,Gravell2020}, which is especially helpful in data-driven applications.

For multiplicative noise the \emph{model-based} setting has already been investigated in \cite{Coppens2019, Coppens2020, Coppens2022, Xing2021}. 
Policy gradient was investigated in \cite{Gravell2021} and policy iteration in \cite{Wang2018, Gravell2022}. The last two 
results however assumed the possibility of exact policy evaluation, which is only possible for known dynamics.

In this work we provide a data-driven scheme for policy iteration applied to linear dynamics with multiplicative noise. 
The approximate policy evaluation step is based on Q-functions as in \cite{Bradtke1994} with the \emph{instrumental variables}
used for \emph{least squares temporal difference (LSTD)} learning in \cite{Bradtke1996}. Our scheme can both operate in 
an \emph{off-policy} setting -- where data-generation happens with some fixed, pre-determined policy -- and in an \emph{on-policy}
setting, where the last policy iterate is applied to the dynamics with some additive noise to enable exploration. 

The remainder of this paper begins with a problem definition in \cref{sec:problem}. Next, in \cref{sec:approximate-policy-evaluation} a novel approximate policy 
evaluation scheme using Q-functions is derived, which is integrated in a policy iteration scheme in \cref{sec:algorithms}.
We also review an existing model-based scheme \citep{Coppens2022} and a policy gradient scheme \citep{Gravell2021}. In \cref{sec:numerical} we 
we then compare the performance and applicability of these three data-driven control schemes. In \cref{sec:conclusion}
we then conclude the paper and suggest further work.

\paragraph*{Notation:~}
Let $\Re$ denote the reals and $\N$ the naturals. For $Z \in \Re^{m \times n}$ let $\trans{Z}$ denote the transpose and 
let $\nrm{Z}_2$ ($\nrm{Z}_F$) be the spectral (Frobenius) norm. When $Z \in \Re^{n \times n}$ let $\lambda(Z) = (\lambda_1, \dots, \lambda_m)$
denote the vector of eigenvalues in descending order of modulus and $\rho(Z) = |\lambda_1(X)|$ the spectral radius. 
For a vector $x \in \Re^d$ let $\nrm{x}_2$ denote the Euclidean norm. Let $\E[\cdot]$ denote the expectation.

We denote by $\sym{d}$ the set of symmetric $d$ by $d$ matrices and by $\sym{d}_{++}$ ($\sym{d}_{+}$) the
positive (semi)definite matrices. Moreover $X - Y \sgt$ ($\sgeq$) $0$ denotes $X - Y \in \sym{d}_{++}$ ($\sym{d}_+$). 
Let $\sd{d} = (d+1)d/2$ denote the degrees of freedom of elements of 
$\sym{d}$ and $\svec(X) \in \Re^{\sd{d}}$ denotes the symmetrized vectorization \citep[\S{}III.A]{Coppens2022} such that 
$\trans{\svec(X)} \svec(Y) = \tr[X Y]$ and $\nrm{\svec(X)}_2 = \nrm{X}_F$,
with $\unsvec \colon \Re^{\sd{d}} \to \Re^d$ its inverse. 

For matrices of conformable size $X, Y$ we use $[X, Y]$ for horizontal concatenation, let $(x, y)$ 
denote vertical concatenation between (column) vectors and let $I_d \in \Re^{d\times d}$ be the identity.
Elements of matrices $X \in \Re^{m \times n}$ (and vectors $x \in \Re^d$) are indexed using $X_{ij}$ ($x_i$) for 
$i \in \{1, \dots, m\}$, $j \in \{1, \dots, n\}$ ($i \in \{1, \dots, d\}$). For a matrix-valued operator $\op{E} \colon \sym{n} \to \sym{m}$ 
we similarly write $\op{E}_{ij} \colon \sym{n} \to \Re$. Moreover $\adj{\op{E}} \colon \sym{m} \to \sym{n}$ is the adjoint, 
such that $\tr[P \op{E}(X)] = \tr[\adj{\op{E}}(P) X]$ for all $X \in \sym{n}$ and $P \in \sym{m}$.

\section{Problem Setup} \label{sec:problem}
Consider linear systems with input- and state-multiplicative noise, given for all $t \in \N$ as:
\begin{equation} \label{eq:dyn}
  x_{t+1} = \left( \sum_{\ell=1}^{n_w} A_i w_{t, \ell} \right) x_t + \left( \sum_{\ell=1}^{n_w} B_i w_{t,\ell} \right) u_t,
\end{equation}
with $x_t \in \Re^{n_x}$ the state, $u_t \in \Re^{n_u}$ the input and 
$w_t \in \Re^{n_w}$ the i.i.d. random disturbance. 

The goal is to solve 
\begin{equation} \label{eq:lqr}
  \minimize_{u_0, u_1, \dots} \quad \E \left[\sum_{t=0}^{\infty} \trans{x_t} Q x_t + \trans{u_t} R u_t \right],
\end{equation}
subject to \eqref{eq:dyn} with $Q \sgt 0$ and $R \sgt 0$ when only given access to trajectories $\{(x_t, u_t)\}_{t=0}^{T}$ of \eqref{eq:dyn}.
The way in which these trajectories are generated (e.g. what policy to use) will be discussed in \cref{sec:data-generation}.

Often it is more convenient to consider the evolution of the second moments $X_t \dfn \E[x_t \trans{x_t}]$, $Z_t \dfn \E[z_t \trans{z_t}]$,
with $z_t = (x_t, u_t)$ the augmented state, when solving \eqref{eq:lqr}. These follow the dynamics 
\begin{equation} \label{eq:dynsm}
  X_{t+1} = \op{E}(Z_t) \dfn \sum_{i,j=1}^{n_w} W_{ij} [A_i, B_i] Z_t \trans{[A_j, B_j]},
\end{equation}
with $W \dfn \E[w_t \trans{w_t}]$ -- independent of $t$ due to the i.i.d. assumption. The LQR problem \eqref{eq:lqr} then becomes 
\begin{equation} \label{eq:lqrsm}
  J_\star(X_0) \dfn \minimize_{Z_0, Z_1, \dots} \quad \sum_{t=0}^{\infty} \tr[H Z_t],
\end{equation}
subject to $Z_t \sgeq 0$ and \eqref{eq:dynsm} starting from $X_0$ and where $H = \diag(Q, R) \sgt 0$. 
We refer to $J(X_0)$ as the value function. 
Equivalence of \eqref{eq:lqr} and \eqref{eq:lqrsm} was shown in \cite[Prop.~IV.6]{Coppens2022}, where the properties 
of $\op{E}$ -- a \emph{completely positive (CP) operator} -- were studied. 

Since we will often deal with partitionings of matrices like $Z \in \sym{n_z}$ we introduce the subscripts:
\begin{equation} \label{eq:partitioning}
  Z = \begin{bmatrix}
    Z_{xx} & \trans{Z_{ux}} \\
    Z_{ux} & Z_{uu}
  \end{bmatrix},
\end{equation}
where $Z_{xx} \in \sym{n_x}$ and analogously for the other terms. 

% The optimal solution is given in terms of a partitioning of the adjoint of $\op{E}$:
% \begin{equation*}
%   \adj{\op{E}}(\cdot) = \begin{bmatrix}
%     [\adj{\op{E}}(\cdot)]_{x, x} & \trans{[\adj{\op{E}}(\cdot)]_{u, x}} \\
%     [\adj{\op{E}}(\cdot)]_{u, x} & \adj{\op{E}}(\cdot)]_{u, u}
%   \end{bmatrix},
% \end{equation*}
% with $[\adj{\op{E}}(\cdot)]_{x, x} = [\adj{\op{E}}(\cdot)]_{1:n_x, 1:n_x}$ and analogously for the others. 
% The optimal policy is then
% \begin{equation*}
%   Z_t = \policy_\star(X_t) \dfn [I; \, K_\star] X_t \trans{[I; \, K_\star]},
% \end{equation*}
% with $K_\star = -(R + [\adj{\op{E}}(P_\star)]_{u, u})^{-1} [\adj{\op{E}}(P_\star)]_{u, x}$
% and $P_\star$ the solution to the Riccati equation $P_\star = \op{R}(P_\star)$ with
% \begin{equation}
%   \begin{aligned}
%   &\op{R}(P_\star) \dfn Q + [\adj{\op{E}}(P_\star)]_{x, x} \\ 
%    &\quad- \trans{[\adj{\op{E}}(P_\star)]_{u, x}} (R + [\adj{\op{E}}(P_\star)]_{x, x})^{-1} [\adj{\op{E}}(P_\star)]_{u, x}.
%   \end{aligned}
% \end{equation}
% The policy $\policy_\star(X_t)$ is implemented in \cref{eq:dyn} as $u_t = K_\star x_t$. 

\section{Approximate Policy Evaluation} \label{sec:approximate-policy-evaluation}
It is well known that the optimal policy solving \cref{eq:lqr} is linear (cf. \cite[Thm.~IV.8]{Coppens2022}). In terms of moments 
such policies look like
\begin{equation} \label{eq:policy}
  Z_t = \policy(X_t) \dfn \trans{[I, \, \trans{K}]} X_t [I, \, \trans{K}],
\end{equation}
We will refer to $K$ and $\policy$ interchangeably as the policy.

As usual, to derive our policy iteration algorithm we need to consider the Bellman operator \citep[\S1.2]{Bertsekas2022}:
\begin{equation} \label{eq:bellman}
  (\op{T}_\policy J)(X) = \tr[\policy(X) H] + J[\op{E}(\policy(X))], \quad \forall X \sgeq 0.
\end{equation}
and $\op{T} J \dfn \min_{\policy} \op{T}_\policy J$, operating on $J \colon \sym{n_x} \to \Re$. 

If the dynamics are stabilizable (i.e. there exists a mean-square stabilizing controller), then
the optimal value of \cref{eq:lqrsm} is given as the fixed-point of $\op{T}$. 
Moreover the optimal controller is then $\argmin_{\policy} \op{T}_{\policy} J_\star$ of \eqref{eq:bellman} with
 $J_\star = \op{T} J_\star$ \citep[Thm.~IV.8]{Coppens2022}.

In our case we know that the value function is linear 
in the second moment\footnote{It is thus quadratic in the state}, i.e. $J(X) = \tr[XP]$ for some $P \sgt 0$. Plugging 
in this parametrization and \cref{eq:policy} and using the definition of the adjoint gives:
\begin{equation*}
  (\op{T} J)(X) = \min_{K}\, \tr[X [I, \, \trans{K}](H + \adj{\op{E}}(P))\trans{[I, \, \trans{K}]}].
\end{equation*}
The minimizer, stated using the subscripts in \cref{eq:partitioning}, is given as $K_{\star} = -(R + \adj{\op{E}}_{uu}(P))^{-1} \adj{\op{E}_{ux}}(P)$.
Now $J = \op{T}J$ can be written in terms of $P$ as 
\begin{align}
  &P = Q + \adj{\op{E}}_{xx}(P) \label{eq:riccati}  - \trans{\adj{\op{E}}_{ux}(P)} (R + \adj{\op{E}_{uu}}(P))^{-1} \adj{\op{E}_{ux}}(P).
\end{align}
This generalized Riccati equation is derived in detail and solved using a SDP in \cite[Thm.~IV.10]{Coppens2022}.

We intend to find $K_\star$ through \emph{policy iteration}, which repeatedly updates $\pi$ by finding a $\pi_+$ such that 
\begin{equation} \label{eq:policy-iteration}
  \pi_{+} = \argmin_{\pi'} T_{\pi'} J_{\pi}, \text{ with } J_{\pi} = T_{\pi} J_{\pi}. 
\end{equation}
Finding $J_{\pi}$ is called \emph{policy evaluation} and requires solving a Lyapunov equation. 
Again using a parametrization $J_\policy(X) = \tr[P_\policy X]$, we write $J_\pi = T_\pi J_\pi$ as:
\begin{align}
  \tr[X P_\policy] &= \tr[\policy(X) H] + \tr[P_\policy \op{E}(\policy(X))] \label{eq:lyapunov}\\
           &= \tr[X \adj{\policy} \left( H + \adj{\op{E}}(P_\policy) \right)], \nonumber
\end{align}
for all $X$, where we used the definition of the adjoint.
Thus we recover the following Lyapunov equation\footnote{
For a proof of invertibility when $\pi$ is mean-square stabilizing see \cite[Prop.~III.8]{Coppens2022}.}:
\begin{align} 
 P_\policy &= \adj{\policy} \left( H + \adj{\op{E}}(P_\policy) \right) \nonumber \\
  \Leftrightarrow \, P_\policy &= (I - \adj{\pi}(\adj{\op{E}}))^{-1}(\adj{\pi}(H)). \label{eq:lyap-solve}
\end{align}
Several schemes exist to solve such equations using data (cf. \cite[\S6]{Bertsekas2012V2} and \cite{Bradtke1996} for MDPs)
However, even after finding $J_\pi$, computing the minimizer in \cref{eq:policy-iteration} requires the dynamics. 
Hence we use Q-functions.

\subsection{Q-functions} \label{sec:q-functions}
A Q-function in this setting maps $Z$ to some real number. We consider the Bellman operator on Q-functions:
\begin{equation*} 
  (\op{F}_\pi \op{Q})(Z) = \tr[ZH] + \op{Q}(\pi(\op{E}(Z)))
\end{equation*}
and $\op{F}\op{Q} \dfn \min_\pi \op{F}_\pi \op{Q}$. For some policy $\pi$ we define $\op{Q}_\pi$ as 
\begin{equation} \label{eq:q-function}
  \op{Q}_\policy(Z) \dfn \tr[Z H] + J_\policy[\op{E}(Z)],
\end{equation}
where $J_\policy$ solves $J_\policy = \op{T}_\policy J_\policy$ (cf. \cref{eq:lyapunov}).
We can show that this choice of $\op{Q}_\policy$ is a fixed-point of $\op{F}_\pi$:
\begin{align*}
  (\op{F}_{\pi} \op{Q}_{\pi})(Z) &= \tr[ZH] + \tr[\pi(\op{E}(Z)) H] + J_\policy[\op{E}(\pi(\op{E}(Z)))] \\
                       &= \tr[ZH] + (\op{T}_\pi J_\pi)[\op{E}(Z)] \\
                       &= \tr[ZH] + J_{\pi}[\op{E}(Z)] = \op{Q}_\pi(Z),
\end{align*}
where we used the definition of $\op{T}_\pi$ for the second equality and $J_\policy = \op{T}_\policy J_\policy$ 
for the third. 

Given such a $\op{Q}_\pi$ we implement policy iteration as:
\begin{equation*}
  \pi_{+}(X) = \argmin_{\pi'} \, \op{Q}_{\pi}(\pi'(X)),
\end{equation*}
which is equivalent to $\pi_{+} = \argmin_{\pi'} \op{F}_{\pi'} \op{Q}_{\pi}$ (and thus \cref{eq:policy-iteration}) yet requires no knowledge about the dynamics.
Using a linear parametrization $\op{Q}_\pi = \tr[\Theta_\pi Z]$ we get
\begin{equation*}
  K_{+} = -(\Theta_{\policy})_{uu}^{-1} (\Theta_{\policy})_{ux}.
\end{equation*} 
Comparing with $K_\star$ from the previous section further confirms that the optimal gain can be recovered 
by selecting the correct value for $\Theta_\policy$.

% Without going into much detail on convergence we examine the fixed-point of this procedure. 
% If we plug in $J_\policy(X) = \tr[P_\policy X]$ into \eqref{eq:q-function} and use adjoints we get:
% \begin{equation} \label{eq:q-function-param}
%   \op{Q}_\policy(Z) = \tr[\Theta_\policy Z] = \tr[(H + \adj{\op{E}}(P_\policy)) Z],
% \end{equation}
% where we define $\Theta = H + \adj{\op{E}}(P_\policy) \sgeq 0$ and partition it as \eqref{eq:partitioning}.
% We can show (similarly to the proof of \cite[Prop.~B.3]{Coppens2022})
% that $K_{k+1}$ -- the value of $K$ for $\policy_{k+1}$ as in \cref{eq:policy} -- is 
% \begin{equation} \label{eq:policy-optimal}
%   K_{k+1} = -[\Theta_{\policy_k}]_{u, u}^{-1} [\Theta_{\policy_k}]_{u, x}.
% \end{equation} 
% If we use $P_\star$ (i.e. the solution of \cref{eq:riccati}) in \cref{eq:q-function-param} we recover 
% $K_\star$ which is then a fixed-point. 

\subsection{Least-Squares} \label{sec:ls}
In our data-driven setting, solving \eqref{eq:lyapunov} to perform the policy evaluation is challenging,
since we need to evaluate $\op{E}$. In \cite{Wang2018} it was assumed that it could be evaluated exactly, which requires knowing the true dynamics.  
Instead, similarly to \cite{Bradtke1994}, we opt to estimate $Q_\pi$ directly from data. We begin by discussing 
constraints on the data and then derive a model equation used for a least-squares estimator. 

We assume access to samples of subsequent moments $\{(Z_i, Z_{i+})\}_{i=1}^{N}$,
where $Z_{i+} = \pi(X_{i+}) = \pi(\op{E}_i(Z_i))$ with
\begin{equation} \label{eq:def-et}
  \op{E}_i(Z_i) \dfn \sum_{j,\ell=1}^{n_w} \left(w_{i, j} w_{i,\ell}\right) [A_j, B_j] Z_i \trans{[A_\ell, B_\ell]}.
\end{equation}
Such data can be generated by taking some $z_i$ and evaluating $x_{i+}$ with \eqref{eq:dyn}.
Then let $Z_i = z_i \trans{z}_i$ and $Z_{i+} = \policy(x_{i+} \trans{x_{i+}})$. We give a detailed data-generation
procedure in \cref{sec:data-generation}.

We want to find $\op{Q}_\pi$ such that $\op{Q}_{\pi}(Z) = (\op{F}_\pi \op{Q}_{\pi})(Z)$ for all $Z \in \psd{n_z}$. 
To do so we sample $\op{F}_{\pi}\op{Q}_{\pi}$:
\begin{equation}
  (\op{F}_{\pi}^i \op{Q}_\pi)(Z_i) \dfn \tr[Z_i H] + \op{Q}_\pi(Z_{i+}),
\end{equation}
which we can use to construct a model equation:
\begin{equation*}
  \op{Q}_\pi(Z_i) = (\op{F}_{\pi}^i \op{Q}_{\pi})(Z_i) + \left((\op{F}_{\pi} \op{Q}_{\pi})(Z_i) - (\op{F}_{\pi}^i \op{Q}_{\pi})(Z_i)\right).
\end{equation*}
Using the parametrization $\op{Q}_{\pi}(Z) = \tr[\Theta_\pi Z]$ results in:
\begin{align}
  \tr[\Theta_\policy Z_i] &=  \tr[Z_i H] + \tr[\Theta_\policy Z_{i+}] \nonumber\\
                &\qquad + \tr[\Theta_\policy \policy(\op{E}(Z_i) - \op{E}_i(Z_i))], \label{eq:model-equation}
\end{align}
which is linear in the parameter $\Theta_\policy$ with zero-mean error. Note that \eqref{eq:model-equation}
involves a \emph{temporal differences} as in \cite{Bradtke1996}.

We reframe the model using the symmetric 
vectorization operator $\svec \colon \sym{n} \mapsto \Re^{\sd{n}}$ with $\sd{n} \dfn n(n+1)/2$ (cf. \cite[\S{}III.A]{Coppens2022}), which satisfies $\tr[X Y] = \trans{\svec{(X)}} \svec{(Y)}$. 
Its inverse is denoted as $\unsvec$. 
We rewrite \eqref{eq:model-equation} as:
\begin{equation}\label{eq:model-equation-vec}
  b_i = \trans{(a_i + e_i)} \theta_\policy,
\end{equation}
with $\theta_\policy = \svec(\Theta_\policy)$, $b_i = \tr[H Z_i]$, $a_i = \svec(Z_i - Z_{i+})$
and $e_i = \svec(Z_{i+} - \policy(\op{E}(Z_{i})))$. This is known as the error-in-variables setting.
The challenge is that both $a_i$ and $e_i$ linearly depend on  $w_i \trans{w_i}$ through $Z_{i+} = \policy(\op{E}_i(Z_i))$. 
Hence $\E[a_i \trans{e}_i] \neq 0$, which makes a classical least-squares estimate inconsistent. 

This is the classical motivation for \emph{instrumental variables (IVs)} \citep[\S3.3]{Young2011}. 
We use $g_i \dfn \svec(Z_i)$ as IVs similarly to \cite{Bradtke1996}. This choice is 
independent of the error $e_i$ (i.e. $\E[g_i \trans{e_i}] = 0$), 
yet is correlated with  $a_i + e_i$ as desired of an IV. 
Intuitively, $g_i$ can be viewed as its best estimate without any additional information on 
the dynamics $\op{E}$. We further motivate the choice by linking it to a projected Bellman equation in \cref{app:ekf}. 

The estimate using the IVs is
\begin{equation*}
  \hat{\theta}_\pi = \left(\sum_{i=1}^{N} g_i \trans{a_i}\right)^{-1} \left(\sum_{i=1}^{N} g_i b_i\right). 
\end{equation*}
Alternatively we can use the recursive algorithm:
\begin{equation} \label{eq:recursive-ivs}
  \begin{aligned}
  \hat{\theta}_{\policy, i} &= \hat{\theta}_{\policy, i-1} + L_i [b_i - \trans{a}_i \hat{\theta}_{\policy, i-1}] \\
  L_i &= S_{\policy, i-1} g_i / (1 + \trans{a_i} S_{\policy, i-1} g_i) \\
  S_{\policy, i} &= S_{\policy, i-1} - (L_i \trans{a_i}) S_{\policy, i-1}
  \end{aligned}
\end{equation}
for $i = 1, \dots, N$ and where the initialization $\hat{\theta}_{\policy, 0} \in \Re^{\sd{n_z}}$ and $S_{\policy, 0} \in \Re^{\sd{n_z} \times \sd{n_z}}$ can be interpreted as 
an initial guess for $\theta_\policy$ and the confidence in that guess respectively. 

\section{Learning Control Schemes} \label{sec:algorithms}
We provide an overview of learning control schemes for the dynamics \cref{eq:dyn}. 
To enable consistent comparisons we first describe a general scheme for gathering data. 
Then we describe our new policy iteration algorithm. Next we briefly summarize the system identification procedure 
of \cite{Coppens2022} and the policy gradient scheme of \cite{Gravell2021}. Sample complexity guarantees 
are reported whenever they exist.

\subsection{Data-generation} \label{sec:data-generation}
We describe the data generation process. 
To enable comparison of the presented algorithms we view them 
as all iteratively updating the policy, where $M$ trajectories (or rollouts) of length $T$ for a total of $N = MT$ new data points 
satisfying \eqref{eq:def-et} are used each iteration. Specifically we simulate in 
parallel for $j = 1, \dots, M$ and $k = 0, \dots, T$:  
\begin{equation*}
  x_{k+1}^{j} = \left( \sum_{\ell=1}^{n_w}A_i w_{k, \ell}^j \right) x_k^{j} +  \left( \sum_{\ell=1}^{n_w} B_i w_{k, \ell}^j \right) u_k^{j},
\end{equation*}
and let $z_{k}^{j} = (x_{k}^{j}, u_{k}^{j})$ for 
\begin{equation} \label{eq:pie}
  u_k^{j} = (K + U^{j}) x_k^{j} + \nu_k^{j}.
\end{equation}
Here $U^{j}$ and $\nu_k^{j}$ are uniformly distributed over $\{U \colon \nrm{U}_F = r_U\}$ and 
$\{\nu \colon \nrm{\nu}_2 \leq r_\nu\}$ respectively. The controller perturbation $U^{j}$ is sampled at 
the start of a trajectory and $\nu_k^{j}$ at each time step. Moreover $U^j$ is distributed over a sphere, instead of a ball,
since in policy gradient, it is used to estimate the gradient of the infinite horizon cost through finite differences as explained in Alg.~\ref{alg:pg}
later. We sample $\nu^j_k$ from a ball to avoid over-excitation of the system, while keeping the trajectories informative enough. 

\begin{figure}
  \usetikzlibrary{hobby,decorations.pathreplacing,arrows.meta}

\tikzset{% adapted from hobby_doc.tex
  show curve controls/.style={
    decoration={
      show path construction,
      curveto code={
        \draw [blue, -{Circle[black,open]}] (\tikzinputsegmentfirst) -- (\tikzinputsegmentsupporta) ;
        \draw [blue, {Circle[black,open]}-] (\tikzinputsegmentsupportb) -- (\tikzinputsegmentlast) ;
      }
    },decorate
  },
}

\begin{tikzpicture}
    \begin{axis}[
        xmin=0, xmax=12,
        ymin=-2.5, ymax=2.5,
        axis x line=bottom,
        axis y line=none,
        height=0.4\columnwidth,
        width=\columnwidth,
        xtick={2.0, 6.0, 10.0},
        xticklabels={iteration $0$, iteration $1$, iteration $2$},
    ]

    % first epoch
    \draw [use Hobby shortcut] 
        ([out angle=20, in angle=180]0.0, 0.4) .. 
        (1.0, 1.0) .. 
        (2.0, 1.5) .. 
        (3.0, 1.2) .. 
        (4.0, 2.0);
    
    \node (v1) at (2.0, 1.5) [circle, fill=red, draw=red, inner sep=0.03cm, node contents={}, label={[label distance=-0.15cm]above right:\footnotesize$x_i$}];
    \node (v2) at (3.0, 1.2) [circle, fill=red, draw=red, inner sep=0.03cm, node contents={}, label={[label distance=-0.15cm]below right:\footnotesize$x_{i+}$}];

    \draw [use Hobby shortcut] 
        ([out angle=-20, in angle=200]0.0, 0.0) .. 
        (1.0, 0.0) .. 
        (2.0, 0.5) .. 
        (3.0, 0.5) .. 
        (4.0, 0.2);
    \draw [use Hobby shortcut] 
        ([out angle=-45, in angle=160]0.0, -0.5) .. 
        (1.0, -1.2) .. 
        (2.0, -1.2) .. 
        (3.0, -1) .. 
        (4.0, -1);

    % second epoch
    \draw [use Hobby shortcut] 
        ([out angle=-20, in angle=180]4.0, 0.6) .. 
        (5.0, 0.5) .. 
        (6.0, 1.0) .. 
        (7.0, 1.5) .. 
        (8.0, 1.0);
    \draw [use Hobby shortcut] 
        ([out angle=0, in angle=160]4.0, -0.1) .. 
        (5.0, -0.2) .. 
        (6.0, 0.5) .. 
        (7.0, 0.5) .. 
        (8.0, 0.2);
    \draw [use Hobby shortcut] 
        ([out angle=-20, in angle=180]4.0, -0.5) .. 
        (5.0, -0.6) .. 
        (6.0, -1.0) .. 
        (7.0, -0.5) .. 
        (8.0, -1);

    % third epoch
    \draw [use Hobby shortcut] 
        ([out angle=45, in angle=180]8.0, 0.4) .. 
        (9.0, 1.0) .. 
        (10.0, 1.0) .. 
        (11.0, 1.5) .. 
        (12.0, 1.0);
    \draw [use Hobby shortcut] 
        ([out angle=20, in angle=160]8.0, -0.2) .. 
        (9.0, 0.5) .. 
        (10.0, 0.2) .. 
        (11.0, 0.5) .. 
        (12.0, 0.7);
    \draw [use Hobby shortcut] 
        ([out angle=-30, in angle=180]8.0, -0.6) .. 
        (9.0, -0.8) .. 
        (10.0, -1.2) .. 
        (11.0, -0.7) .. 
        (12.0, -1.1);
    
    \draw [dashed] (4.0, -2.5) -- (4.0, 2.5);
    \draw [dashed] (8.0, -2.5) -- (8.0, 2.5);

    \end{axis}

\end{tikzpicture}
  \vspace{-0.2cm}
  \centering
  \caption{Rollouts used for data-generation.}\label{fig:trajectories}
\end{figure}
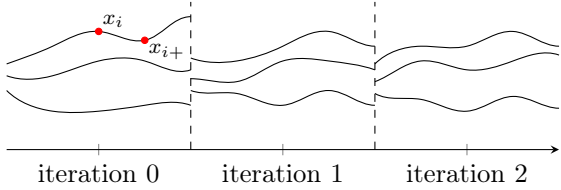

We depict the procedure in \cref{fig:trajectories}.
Since only states and inputs are gathered, the functional form of the dynamics is not required. Only 
a simulator or experiments are required. 

Each trajectory has to be initialized at the start of a new controller iteration. 
We can either continue trajectories from the previous iteration (henceforth referred to as \emph{continuous mode}), or 
we sample $\{x_0^{j}\}_{j=1}^{M}$ uniformly from $\{x \colon \nrm{x}_2 \leq r_x\}$.

We estimate the moments $X_{k}^{j}$ and $Z_{k}^{j}$ for each (augmented) state vector by taking an outer product.
When the specific trajectory does not affect the identification procedure (as in policy iteration and system identification)
we write $(Z_i, X_{i+})$ to index the pairs $(Z_{k}^{j}, X_{k+1}^{j})$, for $k$ and $j$ varying within their range. 
These pairs all satisfy \cref{eq:def-et}.

% \begin{figure}
%   \input{assets/trajectories.tex}
%   \vspace{-0.5cm}
%   \centering
%   \caption{Rollouts used for data-generation. The output for each iteration is $\{Z_i, X_{i+}\}_{i=1}^{N}$
%   with $X_i = x_i \trans{x_i}$, $X_{i+} = x_{i+} \trans{x_{i+}}$ and $Z_{i} = \pi_e(X_i)$.}\label{fig:trajectories}
% \end{figure}

% At each time step we let $Z^{j}_k = \pi_e(X^{j}_k)$ for 
% some exploration policy $\pi_e$, which for fixed state measurements $x_k^{j}$ is implemented as:
% \begin{equation} \label{eq:pie}
%   u_k^{j} = (K + U^{j}) x_k^{j} + \nu_k^{j},
% \end{equation}
% and $Z^{j}_k = \pi_e(X_i) = (x_k^{j}, u_k^{j}) \trans{(x_k^{j}, u_k^{j})}$. 

% Each trajectory in \cref{fig:trajectories} has to be initialized at the start of an iteration. 
% We can either continue trajectories from the previous iteration (henceforth referred to as \emph{continuous mode}), or 
% we sample $x_0$ from some initial state distribution uniformly distributed over $\{x \colon \nrm{x}_2 \leq r_x\}$.

\subsection{Policy Iteration (PI)}

Based on the discussion in the previous section we can now state the full approximate policy iteration algorithm. 
\emph{Generate pair} refers to the rollouts described in the previous section. 

\begin{algorithm2e}[ht!]
  \SetAlgoLined
  \KwData{Initial guess $\hat{\theta}_{\pi_0}$, number of samples $N$ and confidence parameter $\beta_0$.}
  Initialize $S_{\pi_{0}} = \beta_0 I$.\\
  \For{$k=1$ \KwTo $\infty$} 
  {
    $K_{k} = -\Theta_{uu}^{-1} \Theta_{ux}$, with $\Theta = \unsvec(\hat{\theta}_{\pi_{k-1}})$.\\
    Let $\hat{\theta}_{\pi_k, 0} = \hat{\theta}_{\pi_{k-1}}$ and ${S}_{\pi_{k}, 0} = S_{\pi_{k-1}}$. \\
    \For{$i=1$ \KwTo $N$}
    {
      Generate pair $(Z_i, X_{i+})$ and let $Z_{i+} = \pi_{k}(X_{i+})$.\\
      Update $\hat{\theta}_{\pi_{k}, i}$, $S_{\pi_{k}, i}$ using \cref{eq:recursive-ivs}. 
    }
    Set $\hat{\theta}_{\pi_{k}} = \hat{\theta}_{\pi_k, N}$ and $S_{\pi_{k}} = S_{\pi_k, N}$.
  }
  \caption{Approximate policy iteration.} \label{alg:pi}
\end{algorithm2e}

Note how we use the confidence at the end of the previous policy iteration $S_{\pi_k, N}$ 
to initialize the next. This on a high-level corresponds to keeping a summary of the data from past 
policies and is contrasted by \cite{Bradtke1994}, where $S_{\pi_k, 0}$ in \cref{eq:recursive-ivs} is reset to $\beta_0 I$ 
at the start of each policy iteration. Our approach enables the scheme to keep improving 
its estimate of the Q-function as is confirmed in the numerical experiments.
To theoretical back this decision we provide an alternative interpretation as an Extended Kalman Filter 
\citep[\S2.4.3]{Bertsekas2016} applied to a projected Bellman equation in \cref{app:ekf}

\subsection{System Identification (SI)}
The data generated for the policy iteration scheme in the previous section can also be used to identify the dynamics.
We have the following model equation for $i=1, \dots, N$:
\begin{equation*}
  X_{i+} = \op{E}(Z_i) + (\op{E}_i(Z_i) - \op{E}(Z_i)),
\end{equation*}
with $\op{E}_i$ as in \cref{eq:def-et}. The model equation is linear in the matrix representing the linear map $\op{E}$ and has a zero-mean error. Hence we can estimate $\op{E}$
using least-squares as in \cite{Coppens2022}. The resulting least-squares problem has $N \sd{n_x}$ equations and $\sd{n_z}\sd{n_x}$ parameters
describing the matrix of the linear map $\op{E} \colon \sym{n_z} \to \sym{n_x}$. 
The advantage of this scheme is that prior knowledge about the dynamics, like the mode matrices $A_i$ and $B_i$, can be included 
to reduce the problem complexity. In this case we estimate $\E[w \trans{w}]$ instead of $\op{E}$.
Given the dynamics, the optimal controller can be computed by solving a semi-definite program \citep[Thm.~IV.10]{Coppens2022}. 
In the experiments below we do this once per iteration, where the least-square estimate is updated recursively using \cite[Prop.~1.5.2, Eq. 1.118]{Bertsekas2016}. 

Theoretical guarantees are provided in \cite[Thm~VI.4, Rem.~VI.5]{Coppens2022}, which states\footnote{The theoretical guarantees only hold 
for trajectories of length one or when only the final transition of each trajectory is used. The rate was verified empirically for more 
general settings.} that for a
sub-optimality $\epsilon$ the required number of samples is of order $1/\epsilon$. 

\subsection{Policy Gradient (PG)}
We summarize the model-free policy gradient scheme described in \cite{Gravell2021}. 
Given some $\eta$ and initial guess $K_0$ we can compute the optimal controller directly via natural gradient descent:
\begin{equation} \label{eq:gradient-descent}
  K_{k+1} = K_{k} - \eta \nabla \widehat{J(K)} \Sigma_K,
\end{equation}
where the gradient and $\Sigma_K \dfn \sum_{t=0}^{\infty} X_t$ with $X_t$ the closed-loop trajectory with gain $K$
is computed in data-driven fashion using the following algorithm:

\begin{algorithm2e}[ht!]
  \SetAlgoLined
  \KwData{Gain matrix $K$, number of rollouts $M$, rollout length $T$ and exploration radius $r_U$.}
  Generate $\{Z^{j}_t\}$ for $j = 1, \dots M$, $k = 1, \dots, T$ with policy \eqref{eq:pie} for $r_\nu = 0$, $r_U$ and $K$ as provided.\\
  Let $\widehat{J}_j = \sum_{t=0}^{T} \tr[Z_t^{j} H]$ and $\widehat{\Sigma}^j = \sum_{t=0}^{T} (Z_t^j)_{xx}$. \\
  \Return $\nabla \widehat{J(K)} = \sum_{j=1}^{M}\frac{n_x n_u}{Mr^2} \widehat{J}_j U_j$, $\widehat{\Sigma}_K = \frac{1}{M} \sum_{j=1}^M \widehat{\Sigma}^j$. 
  \caption{Gradient evaluation.} \label{alg:pg}
\end{algorithm2e}

Note that this algorithm requires \emph{on-policy} exploration. Moreover multiple rollouts should be generated 
to get a suitable estimate of $\widehat{J}_j$ in each iteration. So policy gradient restricts data-generation more than the 
other methods. 

Theoretical guarantees for classical gradient (i.e. using $\Sigma_K = I$) 
are available in \cite[Thm.~5.1]{Gravell2021}. Given some desired accuracy $\epsilon$,
one should select $M = \mathcal{O}(1/\epsilon^2)$, $T = \mathcal{O}(1/\epsilon^2)$ and $r_U = \mathcal{O}(\epsilon)$
in the gradient evaluation step. The gradient descent scheme will then converge at a linear rate to 
some $K$ such that the sub-optimality is bounded by $\epsilon$. Hence the total sample complexity is of order $1/\epsilon^4$.
In later experiments we illustrate a case with $1/\epsilon$ sample complexity. The stated rates are however 
asymptotic and for classical gradient descent. So the observed gap could disappear for larger sample counts or 
is caused by natural gradient descent outperforming classical gradient descent. Further analysis would thus be interesting. 

\section{Numerical Experiments} \label{sec:numerical}

We will evaluate the presented methods on a dynamical system with multiplicative noise as in \cref{eq:dyn} and modes 
\begin{align*}
  &A_1 = \begin{bmatrix}
      0.43 & 0.71 \\ -1.13 & 0.43
  \end{bmatrix}, \, B_1 = \begin{bmatrix}
    0.36 \\ 0.71
  \end{bmatrix} \, A_2 = \begin{bmatrix}
      0.57 & -0.01 \\ 1.13 & -0.01
  \end{bmatrix}, \\&\quad  B_2=0, \, A_3 = 0,\, B_3 = B_1, \, \text{and } w_t = (1, v_t),
\end{align*}
with $v_t$ sampled uniformly from 
an ellipsoid such that $\E[v_t \trans{v}_t] = \diag(0.2, 0.5)$. The dynamics were motivated by 
a control problem of the pitch of a satellite (cf. \cite[\S1.9.6]{Damm2003}), which 
was discretized using a trapezoid method \citep{Schurz1999}.

We begin with a setup that works best for PI and compare it to SI. Next we compare PG, PI and SI when 
data is generated as prescribed by Alg.~\ref{alg:pg}. Afterwards PI is evaluated in an on-policy setting
and finally we examine the effect of the tuning parameter $\beta_0$ on PI. 

\subsection{Off-policy PI}
First we investigate the effect of the number of rollouts on 
the performance of our policy iteration scheme. Following the notation of \cref{sec:data-generation}, take the additive 
noise radius $r_\nu = 0.1$ and no controller perturbation (i.e. $r_U = 0$). Each 
policy update uses $M = 30$ trajectories of length $T = 100$ are generated,starting from states distributed with $r_x = 1$. 

A fixed, stabilizing control gain $K_0 = [0.5\,-0.75]$
is used to generate the samples. When data is generated using an unstable controller, PI and PG would fail to converge. 
Both schemes would also encounter numerical issues. These can be reduced 
by using more rollouts with a reduced horizon. Further discussion on the effect of stability on PI -- specifically due to 
unstable controllers being generated by the algorithm -- is provided at the end of this section. 

We initialize PI with 
\begin{equation*}
  \hat{\theta}_{\pi_0} = (10, 0, -2.8284, 4, 4.2426, 4), \text{ and } \beta_0 = 2\,000,
\end{equation*}
unless mentioned otherwise. The optimal controller for this initial Q-function equals $K_0$ and is thus stabilizing. 
Note that SI requires no initial guess, which can be considered an advantage.

\begin{figure}[ht!]
  \centering
  \def\svgwidth{\columnwidth}
  \input{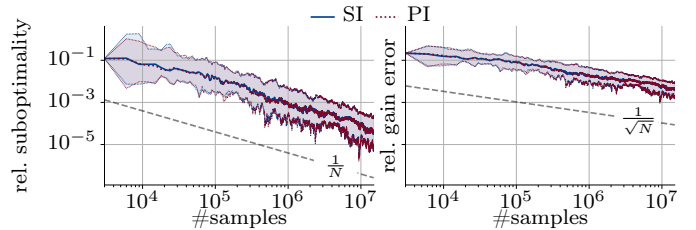}
  \caption{Suboptimality and controller error for off-policy SI and PI with additive exploration noise and $M = 30$.
  The colored area depicts a $80\%$ two-sided empirical confidence interval.} \label{fig:add}
\end{figure}

We terminate after $5000$ iterations and repeat 
the whole process $25$ times to evaluate variances. \cref{fig:add} depicts the comparison between SI and PI. 
For each control gain $K$ and associated $\pi$ we solve the Lyapunov equation \cref{eq:lyap-solve} for the true $\op{E}$
to find value function $J_\pi$. We then compare $J_\pi(I)$ to the optimal $J_\star(I)$ from \cref{eq:lqrsm}. 
This is equivalent to comparing the trace of the associated hessians or the expected infinite closed-loop horizon cost 
starting from a random $x_0$ with $\E[x_0 \trans{x_0}] = I$. The top plots depict the relative error. The bottom 
plots show $\nrm{K - K_\star}_2 / \nrm{K_\star}_2$, i.e. the relative spectral norm error with the optimal gain.

Observe how the suboptimality acts like $1/N$ with $N$ the number of samples, while the controller error acts like $1/\sqrt{N}$
for both PI and SI (the dashed lines follow these rates exactly). 
% Moreover, comparing the left and right plots, 
% we see that increasing the number of rollouts by two orders of magnitude shifts the error down by two orders for suboptimality and one for controller error. 
% So overall the depencency on the number of samples $N$ processed is $1/N$ and $1/\sqrt{N}$ for the suboptimality and controller error respectively.
This corresponds to the rates predicted in \cite[Thm.~VI.4]{Coppens2022}.

For future experiments we will depict the suboptimality only. The
control error -- which decreases at $1/\sqrt{N}$ for all experiments whenever the suboptimality decreases at $1/N$ -- is omitted hereafter.

\subsection{Evaluation of PG}
\begin{figure}[ht!]
  \centering
  \def\svgwidth{0.9\columnwidth}
  \input{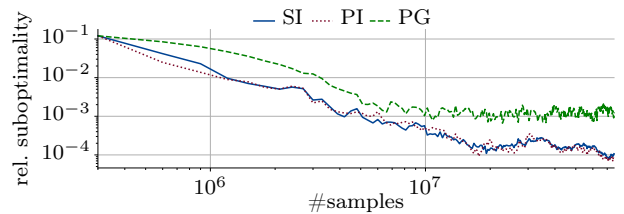}
  \caption{Median relative suboptimality SI, PI and PG.} \label{fig:pg-compare}
\end{figure}

Next we consider natural PG. We tune the data-generation such that PG performs well and then pass the data to SI and PI 
to compare. Here $r_\nu = 0$ and $r_U = 0.15$ produced good results for PG. 
As before we take $M = 3000$ and $T = 100$. We initialize PG with $K_0 = [0.5\,-0.75]$ and use the PG gain in later 
iterations as described in Alg.~\ref{alg:pg}. The step size is set to $7.5\cdot 10^{-3}$. We initialize PI 
with $\beta_0 = 5.0$. This change compensates for the absence of additive noise, which makes the data less informative.
It corresponds to more trust in the initial guess to avoid generating unstable controllers, which causes PI to diverge. 
Another experiment below expands on this intuition. We continue for $250$ iterations. 
The rest of the procedure is identical to before and the result is depicted in \cref{fig:pg-compare}.

\begin{figure}[ht!]
  \centering
  \def\svgwidth{0.9\columnwidth}
  \input{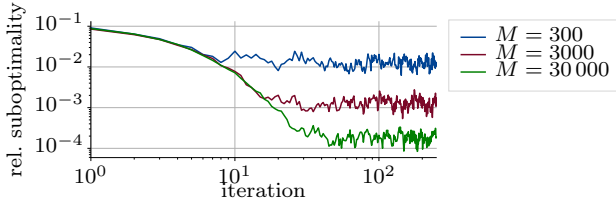}
  \caption{Evaluation of sample complexity of PG. Plots show median relative suboptimality.} \label{fig:pg-samples}
\end{figure}

Note how SI and PI act similarly to before, while PG stagnates after a certain number of iterations. 
This is also predicted by the theory in \cite{Gravell2021}. We can evaluate the sample complexity 
by repeating the same experiment, but with $M = 300$, $3000$ and $30\,000$. We plot the number of iterations on the horizontal axis to 
ease visual comparison, keeping in mind that the number of samples observed per iteration differs. The result in \cref{fig:pg-samples} 
indicates that the suboptimality improves with $1/M$. So overall the rate is the same as SI and PI, but it 
manifests in terms of the number of rollouts per gradient evaluation, instead of the cumulative sample count. Moreover, while SI and PI put no requirements on data generation,
PG is more restrictive.

\subsection{On-policy PI}
\begin{figure}[ht!]
  \centering
  \def\svgwidth{0.9\columnwidth}
  \input{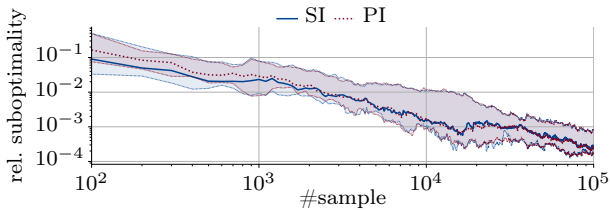}
  \caption{On-Policy PI and SI.
  The colored area depicts a $80\%$ two-sided empirical confidence interval.} \label{fig:on-policy}
\end{figure}

We can run PI in an on-policy setting with one rollout each iteration. Taking $T=100$, $M=1$ 
and $x_0$ distributed with $r_x = 1$, while using the final state of the latest iteration afterwards (c.f. continuous mode in \cref{sec:data-generation}). 
We initialize PI with $\beta_0 = 100$ and run the scheme for $1000$ iterations. The rest of the setup is identical to the first experiment, but now instead of applying $K_0$ 
for every iteration, the latest PI control gain is used instead. The result, manifesting the same sample complexity as before, 
is depicted in \cref{fig:on-policy}. Interestingly the suboptimality is reduced compared to \cref{fig:add}, potentially
indicating that this method of data-generation is more informative.

\subsection{Tuning of PI}
When a small number of samples are provided per iteration in PI, unstable control gains can be generated. 
This can cause the algorithm to diverge.\footnote{See \cite[\S3]{Bertsekas2022}, which proves convergence of PI over 
the class of stable policies and illustrates what happens when an unstable policy is used instead.} 
To examine the effect of tuning we set $T = 10$ and $M = 5$ and run the algorithm $1000$ times for $100$ iterations. 
The remainder of the setup is the same as the first experiment of this section. The data-generation uses the
same stabilizing policy $K_0 = [0.5 \, -0.75]$ for each iteration.

\begin{figure}[ht!]
  \centering
  \def\svgwidth{0.9\columnwidth}
  \input{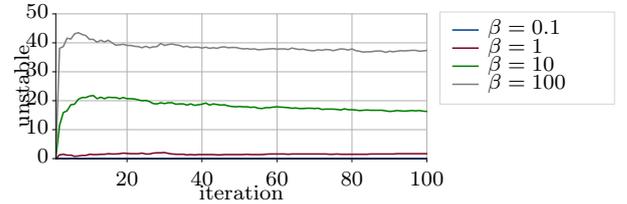}
  \caption{Percentage of unstable policies produced by PI for different tunings.} \label{fig:pi-fail}
\end{figure}

The percentage of unstable policies per iteration is depicted in \cref{fig:pi-fail}. Note how unstable policies are produced 
in the first iterations due to the small amount of data and mostly remain unstable afterwards.
Low values for $\beta_0$ reduces the number of unstable policies since the algorithm trusts the initial guess more. 
For $\beta_0 = 0.1$ no unstable policies occur. Alternatively, performing system identification separately from control 
synthesis does not encounter stability issues and requires no tuning. 

\section{Conclusion} \label{sec:conclusion}
From the experiments in the previous section one can conclude a pay-off between a simple implementation and 
robustness. The policy gradient scheme requires the most effort to tune and requires the most 
specific data, yet is simple to implement. Meanwhile the model-based approach requires solving a large least-squares problem and 
a generalized algebraic Riccati equation -- usually formulated as a SDP. However it requires no tuning and can
handle almost all data-generation schemes. The presented policy iteration scheme meanwhile also solves a smaller least-squares problem,
yet avoids the Riccati equation. However as shown in \cref{fig:pi-fail}, with limited data the scheme can require 
some tuning to work. 

In future work we will investigate adaptive iteration lengths to avoid such tuning issues. After all,
problems often occur in the first iterations. Hence gathering more data there would aid in stabilizing the algorithm. 
The connection with EKF in \cref{app:ekf} will also be exploited further to find stability and convergence guarantees.
Such a theoretical analysis would also aid in improving tuning.

\begin{ack}
  We thank Jean-Louis Carron for performing an initial comparison of learning control schemes for multiplicative
  noise during their master thesis, motivating the developments presented here.
\end{ack}

\printbibliography                  % bib file to produce the bibliography
                                    % with bibtex (preferred)

\appendix

\section{EKF Interpretation} \label{app:ekf}
\paragraph*{Introduction to Extended Kalman Filters:~}
Classically an EKF is applied to a nonlinear least-squares problem:
\begin{equation*}
  \minimize_{\theta} \quad \frac{1}{2} \nrm{g(\theta)}_2^2 = \frac{1}{2} \sum_{i=1}^{N} \nrm{f_i(\theta)}_2^2. 
\end{equation*}
It linearizes $f_i(\theta)$ around recursive estimates of the minimizer $\theta_k$, denoted as 
$\tilde{f}_i(\theta | \theta_k) \dfn f_i(\theta_k) + \trans{\nabla f_i(\theta_k)} (\theta - \theta_k)$. 
These estimates are then computed recursively as \citep[Eq.~2.129]{Bertsekas2012V2}:
\begin{equation}
  \theta_k \in \argmin_{\theta} \, \sum_{\ell=1}^{k} \nrm{\tilde{f}_\ell(\theta \mid \theta_{\ell-1})}_2^2, \quad k = 1,\dots, N.
\end{equation}
Our scheme applies the same idea instead to a projected equation of the form 
\begin{equation} \label{eq:projected-general-equation}
  \hat{\theta} \in \argmin_{\theta} \, \sum_{i=1}^{M} \nrm{\trans{\varphi_i} \theta - h_i(\hat{\theta})}_2^2,
\end{equation}
where we now linearize $h_i$ at each step, producing an algorithm like
\begin{equation} \label{eq:general-ekf}
  \hat{\theta}_k \in \argmin_{\theta} \, \sum_{\ell=1}^{k}  \nrm{\trans{\varphi_i} \theta - \tilde{h}_i(\hat{\theta}_k \mid \theta_{\ell-1})}_2^2.
\end{equation}
Instead of linearizing each time step however we split the data up into batches and linearize at the start of each batch.
Next, this reasoning is applied to our setting.

\paragraph*{EKF Policy Iteration:~}
We reinterpret Alg.~\ref{alg:pi} as an \emph{Extended Kalman Filter} applied to a projected Bellman equation. 
Following the notation of \cref{sec:q-functions}, a policy update looks like:
\begin{equation*}
  \pi_{+} = \argmin_{\pi'} \op{F}_{\pi'} \op{Q}_{\pi}, \text{ with } \op{Q}_{\pi} = \op{F}_{\pi} \op{Q}_{\pi}.
\end{equation*}
Note that $\op{F}_{\pi_+}$ is linear in $\op{Q}$ and $\op{F}_{\pi_+} \op{Q}_{\pi} = \op{F} \op{Q}_{\pi}$. 
Thus
\begin{align*}
  \mathcal{F} \op{Q}_\pi + \mathcal{F}_{\pi_+}(\op{Q} - \op{Q}_\pi) = \mathcal{F}_{\pi_+}\op{Q} \geq \mathcal{F} \op{Q}, \quad \forall \op{Q}.
\end{align*}
So $\mathcal{F}_{\pi_+}$ acts as a supergradient of $\mathcal{F}$ and
 $\mathcal{F}_{\pi_+}(\op{Q})$ is a linearization of $\mathcal{F}(\op{Q})$ at $\op{Q}_\pi$. 

Similar to \cref{sec:data-generation} we consider a linear parametrization of $\op{Q}$
\begin{equation*}
  \op{Q}^\theta(Z) \dfn \trans{\svec}(Z) \theta
\end{equation*}
and introduce the sampled Bellman operator:
\begin{equation*}
  (\mathcal{F}_\pi^i \op{Q})(Z) = \tr[ZH] + \op{Q}(\pi(\op{E}_i(Z))).
\end{equation*}
and let $\mathcal{F}^i \op{Q} = \min_{\pi} \mathcal{F}_{\pi}^i \op{Q}$.

Consider now the projected Bellman equation\footnote{Here $\op{Q}^{\hat{\theta}_{\pi}}$ is the parametrized $\op{Q}_\pi$.} \citep[\S6.3]{Bertsekas2012V2}:
\begin{equation} \label{eq:projected-bellman}
  \hat{\theta}_\pi = \argmin_{\theta} \sum_{i=1}^{N} \nrm{\op{Q}^\theta(Z_i) - (\mathcal{F}^i \op{Q}^{\hat{\theta}_\pi})(Z_i)}_2^2.
\end{equation}
This is our specific case of \cref{eq:projected-general-equation}. We can view Alg.~\ref{alg:pi} as 
solving \cref{eq:projected-bellman} through a scheme reminiscent of an EKF. Specifically we split the data into batches of size $N_k$ 
denoted as $\mathcal{I}_k = (\bar{N}_{k-1} + 1, \bar{N}_{k-1} + 2, \dots, \bar{N}_{k-1} + N_{k})$, with $\bar{N}_k = \sum_{\ell=1}^k N_\ell$
and $\bar{N}_0 = 0$. 

The algorithm then updates $\hat{\theta}_{\pi_k}$ each iteration as (cf. \cref{eq:general-ekf}):
\begin{align*}
  \hat{\theta}_{\pi_k} = \argmin_{\theta} \left\{ \sum_{\ell=1}^k \sum_{i \in \mathcal{I}_\ell} \nrm{\op{Q}^\theta(Z_{i}) - (\mathcal{F}^{i}_{\pi_{\ell}} \op{Q}^{\hat{\theta}_{\pi_k}})(Z_{i})}_2^2 \right\},
\end{align*}
where we linearized the term in the norm of \cref{eq:projected-bellman} by
taking $\pi_{\ell}$ such that $\mathcal{F} \op{Q}^{\hat{\theta}_{\pi_{\ell-1}}} = \mathcal{F}_{\pi_{\ell}}  \op{Q}^{\hat{\theta}_{\pi_{\ell-1}}}$. 
Unlike in classical EKF we linearize only at the start of each batch, not for each sample. This aids in stabilizing
the algorithm.

To make the connection with Alg.~\ref{alg:pi} complete we re-introduce the notation of \cref{sec:ls}.
Specifically $b_i = \tr[H Z_i]$, $g_i = \svec(Z_i)$, $g_{i+} = \svec(\pi_{\ell-1}(X_{i+1}))$ and $a_i = g_i - g_{i+}$. Then we get:
\begin{equation} \label{eq:projected-equation-ls}
  \hat{\theta}_{\pi_k} = \argmin_{\theta} \left\{\sum_{i=1}^{\bar{N}_k} \nrm{\trans{g_i} \theta - \trans{g_{i+}} \hat{\theta}_{\pi_k} - b_i}_2^2 \right\},
\end{equation}
which is equivalent to the following normal equations:
\begin{align*} 
  \sum_{i=1}^{\bar{N}_k} g_i((\trans{g_i} - \trans{g_{i+}}) \hat{\theta}_{\pi_k}- b_i) = \sum_{i=1}^{\bar{N}_k} g_i(\trans{a_i} \hat{\theta}_{\pi_k} - b_i) = 0.
\end{align*}
This corresponds to the normal equation of least-squares with IVs, hence motivating their use in \cref{sec:ls}. We 
can use the recursion \cref{eq:recursive-ivs} to generate solutions iteratively without computing inverses. 

Specifically we add a regularization term \[\nrm{\theta - \hat{\theta}_{\pi_0}}_{S_{\pi_0}^{-1}}^2/2\] to \cref{eq:projected-equation-ls}.
The updates can then be computed recursively by solving a one-stage normal equation each time:
\begin{equation} \label{eq:single-stage-ivs}
  g_i (\trans{a}_i \hat{\theta}_{\pi_{k}} - b_{i}) + S_{\pi_{k}, i-1}^{-1}(\hat{\theta}_{\pi_{k}} - \hat{\theta}_{\pi_{k}, {i-1}}) = 0.
\end{equation}
The solution of which is
\begin{equation} \label{eq:ivs-deriv-0}
  \hat{\theta}_{\pi_k, i} = \left( S_{\pi_k, i-1}^{-1} + g_i \trans{a_i} \right)^{-1} (S_{\pi_k, i-1}^{-1} \hat{\theta}_{\pi_k, i-1} + g_i b_i).
\end{equation}
The first factor is simplified through Woodbury's matrix identity applied to the first factor and introducing $L_{i} = S_{\pi, i-1} g_i / (1+ \trans{a_i} S_{\pi_k, i-1} g_i)$
to
\begin{align}\label{eq:ivs-deriv-1}
  S_{\pi_k, i} &\dfn \left( S_{\pi_k, i-1}^{-1} + g_i \trans{a_i} \right)^{-1} \\
  &= S_{\pi_k, i-1} - (L_i \trans{a_i}) S_{\pi_k, i-1}. \nonumber
\end{align}
Plugging into \cref{eq:ivs-deriv-0} gives 
\begin{equation*}
  \hat{\theta}_{\pi_k, i} = \hat{\theta}_{\pi_k, i-1} + L_i[b_i - \trans{a_i} \hat{\theta}_{\pi_k, i-1}].
\end{equation*}
Also, by combining this with \cref{eq:ivs-deriv-0} and \cref{eq:ivs-deriv-1} we can show that \cref{eq:single-stage-ivs}
is equivalent to 
\begin{equation} \label{eq:regularization-final}
  S_{\pi_k, i}^{-1} (\hat{\theta}_{\pi_k} - \hat{\theta}_{\pi_k, i}) = 0.
\end{equation}
So we can equivalently process the next terms in \cref{eq:projected-equation-ls} by adding the 
left-hand side of \cref{eq:regularization-final} as a regularization 
term. As such we have derived the recursion in \cref{eq:recursive-ivs}.

\end{document}